# Some observations about super Catalan numbers, corresponding orthogonal polynomials, and their q-analogues

Johann Cigler

johann.cigler@univie.ac.at

Fakultät für Mathematik, Universität Wien,

Uni Wien Rossau, Oskar Morgenstern-Platz 1, 1090 Wien

**Abstract**

In this note we give a survey about polynomials whose moments are multiples of super Catalan numbers and explore two different kinds of $q-$ analogues.

**0. Introduction**

The moments of Lucas polynomials and of Chebyshev polynomials of the first kind are (multiples of) central binomial coefficients and the moments of Fibonacci polynomials and of Chebyshev polynomials of the second kind are (multiples of) Catalan numbers. In [6] I presented some generalizations of these results together with some $q-$ analogues. The recent paper [1] by E. Allen and I. Gheorghiciuc has prompted me to look for nice polynomials whose moments are super Catalan numbers or $q-$ analogues of super Catalan numbers. This search led to natural extensions of Fibonacci and Chebyshev polynomials and to two different kinds of $q-$ analogues. Some special cases have already occurred in [4] and [6]. Once the results are known proofs are straightforward verifications and will therefore be omitted. Finally I show how the polynomials with odd degree can be expressed by the polynomials with even degree.

To provide some background information let us start with some well-known facts (cf. e.g. [6]). Consider the variant $l_n(x)$ of the Lucas polynomials defined by

$$l_n(x) = \sum_{k=0}^{\lfloor \frac{n}{2} \rfloor} \binom{n-k}{k} \frac{n}{n-k} (-1)^k x^{n-2k} \tag{0.1}$$

$$l_0(x) = 1.$$

The polynomials $l_n(x)$ satisfy the recurrence relation

$$l_n(x) = x l_{n-1}(x) - \lambda_{n-2} l_{n-2}(x) \tag{0.2}$$

with $\lambda_0 = 2$ and $\lambda_n = 1$ for $n > 0$.

Define a linear functional $\Lambda$ on the polynomials by $\Lambda(l_n(x)) = [n = 0]$.

From the representation

$$x^n = \sum_{k=0}^{\lfloor n/2 \rfloor} \binom{n}{k} l_{n-2k}(x) \qquad (0.3)$$

we deduce that $\Lambda(x^{2n+1}) = 0$ and

$$\Lambda(x^{2n}) = \binom{2n}{n}. \qquad (0.4)$$

An analogous situation will occur several times in the sequel. In all cases we will have $\Lambda(x^{2n+1}) = 0$. The numbers $\Lambda(x^{2n})$ will be called the moments of the linear functional $\Lambda$.

Consider next the variant $f_n(x)$ of the Fibonacci polynomials defined by

$$f_n(x) = \sum_{k=0}^{\lfloor n/2 \rfloor} \binom{n-k}{k} (-1)^k x^{n-2k}. \qquad (0.5)$$

Define a linear functional $\Lambda^*$ on the polynomials by $\Lambda^*(f_n(x)) = [n = 0]$.

From the representation

$$x^n = \sum_{k=0}^{\lfloor n/2 \rfloor} \left(\binom{n}{k} - \binom{n}{k-1}\right) f_{n-2k}(x) \qquad (0.6)$$

we deduce that the moments are the Catalan numbers

$$\Lambda^*(x^{2n}) = \frac{1}{n+1}\binom{2n}{n} = C_n. \qquad (0.7)$$

The polynomials $l_n(x)$ and $f_n(x)$ satisfy $l_n(x) = f_n(x) - f_{n-2}(x)$ for $n \geq 2$.

Closely related are the monic Chebyshev polynomials.

The monic Chebyshev polynomials of the first kind $t_n(x)$ satisfy

$$t_n(x) = x t_{n-1}(x) - \lambda_{n-2} t_{n-2}(x) \tag{0.8}$$

with $\lambda_0 = \dfrac{1}{2}$ and $\lambda_n = \dfrac{1}{4}$ for $n > 0$. Their moments are $\dfrac{1}{2^{2n}}\binom{2n}{n}$.

The monic Chebyshev polynomials of the second kind $u_n(x)$ satisfy

$$u_n(x) = x u_{n-1}(x) - \lambda_{n-2} u_{n-2}(x) \tag{0.9}$$

with $\lambda_n = \dfrac{1}{4}$ for $n \geq 0$. Their moments are $\dfrac{1}{2^{2n}} C_n$.

## 1. Super Catalan numbers and related polynomials

Consider now the super Catalan numbers $S(m,n) = \dfrac{(2m)!(2n)!}{m!n!(m+n)!}$.

These numbers and their $q-$analogues have been studied in [1], [9] and [11].

We look for monic polynomials whose moments are related to $S(m,n)$. Since we want $\Lambda(1) = 1$ we must normalize the sequence $S(m,n)$, i.e. consider the sequence $(\sigma(m,n))_{n \geq 0}$ of the "normalized super Catalan numbers"

$$(\sigma(m,n))_{n \geq 0} = \frac{S(m,n)}{S(m,0)} = \frac{(2n)!m!}{n!(n+m)!}. \tag{1.1}$$

For $m = 0$ we get the central binomial coefficients $(\sigma(0,n))_{n \geq 0} = \dfrac{(2n)!}{n!(n)!} = \binom{2n}{n}$ and for $m = 1$ we see that $(\sigma(1,n))_{n \geq 0} = \dfrac{(2n)!}{n!(n+1)!} = C_n$. These are the only cases where all terms are integers.

The first terms of the sequence $((\sigma(m,n))_{n \geq 0})_{n \geq 0}$ are

$$\left(1, \frac{2}{m+1}, \frac{12}{(m+1)(m+2)}, \frac{120}{(m+1)(m+2)(m+3)}, \frac{1680}{(m+1)(m+2)(m+3)(m+4)}, \ldots\right).$$

For example $\big((\sigma(2,n))_{n\geq 0}\big)_{n\geq 0} = \left(1, \dfrac{2}{3}, 1, 2, \dfrac{14}{3}, 12, 33, \dfrac{286}{3}, \cdots\right).$

The search for polynomials whose moments are $\sigma(m,n)$ led me to the sequence $l_n(x,m,s)$ defined by

$$l_n(x,m,s) = \sum_{k=0}^{\lfloor \frac{n}{2} \rfloor} \frac{n!}{(k)!(n-2k)!} \frac{1}{\prod_{j=1}^{k}(m+n-j)} s^k x^{n-2k}. \tag{1.2}$$

The first terms are

$1,\ x,\ x^2 + \dfrac{2s}{m+1},\ x^3 + \dfrac{6s}{m+2}x,\ x^4 + \dfrac{12s}{m+3}x^2 + \dfrac{12s^2}{(m+2)(m+3)},\ x^5 + \dfrac{20s}{m+4}x^3 + \dfrac{60s^2}{(m+3)(m+4)}x, \cdots.$

By comparing coefficients it is easy to show that these polynomials satisfy the recursion

$$l_n(x,m,s) = x l_{n-1}(x,m,s) + s\lambda_{n-2}(m) l_{n-2}(x,m,s) \tag{1.3}$$

where $\lambda_0(m) = \dfrac{2}{m+1}$ and for $n > 0$

$$\lambda_n(m) = \dfrac{(n+1)(n+2m)}{(n+m)(n+m+1)}. \tag{1.4}$$

Since they satisfy the 3-term recurrence (1.3) they are orthogonal with respect to the linear functional $\Lambda_m$ on the polynomials defined by $\Lambda_m\big(l_n(x,m,s)\big) = [n=0]$.

For $m = 0$ we get $l_n(x,0,-1) = l_n(x)$.

For $m = 1$ we have $l_n(x,1,-1) = f_n(x)$.

Therefore $l_n(x,m,-1)$ is a natural extension of Lucas and Fibonacci polynomials for all $m \in \mathbb{N}$.

If we choose $s = -\dfrac{1}{4}$ we get the monic Chebyshev polynomials $l_n\left(x,0,-\dfrac{1}{4}\right) = t_n(x)$ and $l_n\left(x,1,-\dfrac{1}{4}\right) = u_n(x)$.

Let us note some special values:

The sequence

$$(l_n(1,0,-1))_{n\geq 0} = (1,1,-1,-2,-1,1,2,1,-1,-2,-1,1,\cdots) \quad (1.5)$$

is periodic with period 6 if we replace $l_0(1,0,-1)$ by 2.

The sequence

$$(l_n(1,1,-1))_{n\geq 0} = (1,1,0,-1,-1,0,1,1,0,-1,-1,0,\cdots) \quad (1.6)$$

is also periodic with period 6.

For the sequence $(l_n(1,2,-1)(n+1))_{n\geq 0}$ we get

$$(l_n(1,2,-1)(n+1))_{n\geq 0} = (1,2,1,-2,-4,-2,3,6,3,-4,-8-,4,\cdots), \quad (1.7)$$

i.e. $l_{3n}(1,2,-1)(3n+1) = l_{3n+2}(1,2,-1)(3n+3) = (-1)^n(n+1)$,

and $l_{3n+1}(1,2,-1)(3n+2) = (-1)^n 2(n+1)$.

For the Chebyshev-like polynomials we get

$$l_n\left(1,m,-\frac{1}{4}\right) = \frac{(2m+1)(2m+2)\cdots(2m+n-1)}{(2m+2)(2m+4)\cdots(2m+2n-2)} = \prod_{j=1}^{\lfloor \frac{n}{2}\rfloor} \frac{2m+2j-1}{2\left(m+j+\left\lfloor\frac{n-1}{2}\right\rfloor\right)}. \quad (1.8)$$

The uniquely determined coefficients of the expansion

$$x^n = \sum_{k=0}^{\lfloor \frac{n}{2}\rfloor} a(n,k,m)(-s)^k l_{n-2k}(x,m,s) \quad (1.9)$$

satisfy $a(n+1,k,m) = a(n,k,m) + \lambda_{n+1-2k}(m)a(n,k-1,m)$. This gives by induction

$$a(n,k,m) = \binom{n}{k}\frac{(n-k)!(n+m-2k)!}{(n-2k)!(n+m-k)!} = \frac{n!}{k!(n-2k)!}\frac{(n+m-2k)!}{(n+m-k)!}. \quad (1.10)$$

Thus we have

$$x^n = \sum_{k=0}^{\lfloor \frac{n}{2}\rfloor} \frac{n!}{k!(n-2k)!}\frac{(n+m-2k)!}{(n+m-k)!}(-s)^k l_{n-2k}(x,m,s) \quad (1.11)$$

which implies that the moments are multiples of the normalized super Catalan numbers

$$\Lambda_m\left(x^{2n}\right) = (-s)^n \frac{(2n)!m!}{n!(m+n)!} = (-s)^n \sigma(m,n).$$

Instead of (1.9) we can also consider the coefficients $b(n,k,m)$ in the expansion

$$x^n = \sum_{k=0}^{n} b(n,k,m) l_k(x,m,s). \tag{1.12}$$

They satisfy

$$b(n,k,m) = b(n-1,k-1,m) + s\lambda_k(m) b(n-1,k+1,m) \tag{1.13}$$

since

$$\sum_{k=0}^{n} b(n,k,m) l_k(x,m,s) = x \cdot x^{n-1} = \sum_{k=0}^{n} b(n-1,k,m)\bigl(x l_k(x,m,s)\bigr)$$

$$= \sum_{k=0}^{n} b(n-1,k,m)\bigl(l_{k+1}(x,m,s) + s\lambda_{k-1}(m) l_{k-1}(x,m,s)\bigr)$$

$$= \sum_{k} b(n-1,k-1,m) l_k(x,m,s) + \sum_{k} b(n-1,k+1,m) s\lambda_k(m) l_k(x,m,s).$$

For $m=1$ these are the so-called ballot numbers.

Comparing (1.12) with (1.9) we see that

$$b(2n,2k,m) = a(2n,n-k,m)(-s)^{n-k} = \frac{(2n)!(m+2k)!(-s)^{n-k}}{(n-k)!(2k)!(n+m+k)!} = \binom{2n}{2k}\sigma(m+2k,n-k)(-s)^k.$$
$$\tag{1.14}$$

Let us state some simple well known facts about generalized ballot numbers.

**Lemma 1**

*Let $b(n,k)$ satisfy*

$$b(n,k) = b(n-1,k-1) + r(k) b(n-1,k+1) \tag{1.15}$$

*for $n > 0$ where $b(n,-1) = 0$ and $b(0,k) = [k = 0]$.*

*Then the numbers $b(n,k)$ can be interpreted as the weight of the set of all non-negative lattice paths from $(0,0)$ to $(n,k)$ consisting of up-steps $(i,j) \to (i+1, j+1)$ with weight $1$ and down-steps $(i, j+1) \to (i+1, j)$ with weight $r(j)$.*

*Let $c(\ell,k)$ be the weight of the set of all non-negative lattice paths from some point $(n,k)$ to $(n+\ell,0)$. Then $c(2\ell, 2k) = b(2\ell, 2k) \prod_{j=0}^{2k-1} r(j)$ since for each height $j$ the number of down-steps $(i, j+1) \to (i+1, j)$ is one more than the number of up-steps $(i,j) \to (i+1, j+1)$.*

*Therefore we have*

$$\sum_{k=0}^{\min(n,\ell)} b(2n, 2k) b(2\ell, 2k) \prod_{j=0}^{2k-1} r(j) = b(2n+2\ell, 0). \tag{1.16}$$

Thus (1.16) implies for $s = -1$

$$\sum_{k=0}^{\min(n,\ell)} b(2n, 2k, m) b(2\ell, 2k, m) \prod_{j=0}^{2k-1} \left(\lambda_j(m)\right) = \sigma(m, n+\ell). \tag{1.17}$$

For $m=0$ this reduces to

$$\binom{2n}{n}\binom{2\ell}{\ell} + 2 \sum_{k=1}^{\min(n,\ell)} \binom{2n}{n+k}\binom{2\ell}{\ell+k} = \binom{2n+2\ell}{n+\ell} \tag{1.18}$$

or equivalently to

$$\sum_{k=-n}^{n} \binom{2n}{n-k}\binom{2\ell}{\ell-k} = \binom{2n+2\ell}{n+\ell}. \tag{1.19}$$

For $m=1$ we get

$$\sum_{k=0}^{n} \binom{2n}{n-k}\binom{2\ell}{\ell-k} \frac{(2k+1)^2}{(n+k+1)(\ell+k+1)} = \binom{2n+2\ell}{n+\ell} \frac{1}{n+\ell+1}. \tag{1.20}$$

Let us also state the following result (cf. [4]), which can be verified by induction:

**Lemma 2**

*Under the same assumptions as above*

$$\sum_{k=0}^{n}(-1)^k b(2n,2k)\prod_{j=0}^{k-1} r(2j) = [n=0]. \tag{1.21}$$

By induction we have for $m \geq 0$

$$\prod_{j=0}^{k-1} \lambda_{2j}(m) = \sigma(m+k-1,k). \tag{1.22}$$

For $m=0$ this reduces to $\left(\prod_{j=0}^{k-1} \lambda_{2j}(0)\right)_{k\geq 0} = (1,2,2,2,\cdots)$ and (1.21) to

$$\sum_{k=-n}^{n}(-1)^k \binom{2n}{n-k} = [n=0]. \tag{1.23}$$

For $m=1$ we have $\sigma(k,k)=1$ and therefore

$$\sum_{k=0}^{n}(-1)^k \binom{2n}{n-k}\frac{2k+1}{n+k+1} = [n=0]. \tag{1.24}$$

For $m>0$ formula (1.21) can be simplified to

$$\sum_{k=0}^{n}(-1)^k \binom{n}{k}\frac{(m+2k)(m+k-1)!}{(n+m+k)!} = [n=0]. \tag{1.25}$$

In general there is no simple formula for $\sum_{k=0}^{n} b(2n,2k)\prod_{j=0}^{k-1} r(2j)$, but in our special case we get

$$\sum_{k=0}^{n} b(2n,2k,m)\prod_{j=0}^{k-1} \lambda_{2j}(m) = \frac{2^{2n}(2n-1)!!}{(m+1)(m+3)\cdots(m+2n-1)}. \tag{1.26}$$

For $m>0$ it can be simplified to

$$\sum_{k=0}^{n}\binom{n}{k}\frac{(m+2k)(m+k-1)!}{(n+m+k)!} = \frac{2^n}{(m+1)(m+3)\cdots(m+2n-1)}. \tag{1.27}$$

A simple computer proof can be given by using the Mathematica implementation of Zeilberger's algorithm fastZeil by Peter Paule and Markus Schorn which gives for the left-hand side $S(n)$ the recurrence $2S(n)-(m+1+2n)S(n+1)=0,$ which immediately implies the right-hand side.

**Remark**

For $m>0$ a generating function of the polynomials $l_n(x,m,s)$ is given by

$$\sum_{n\geq 0} l_n(x,m,s)\binom{n+m-1}{m-1}z^n = \frac{1}{(1-xz-sz^2)^m} = \left(\sum_{n\geq 0} l_n(x,1,s)z^n\right)^m. \quad (1.28)$$

To verify this identity differentiate it with respect to $z$ and multiply with $z(1-z-sz^2)$. This gives

$$(1-xz-sz^2)\sum_{n\geq 0} l_n(x,m,s)\binom{n+m-1}{m-1}nz^n = m\frac{xz+2sz^2}{(1-xz-sz^2)^m} = m(xz+2sz^2)\sum_{n\geq 0} l_n(x,m,s)\binom{n+m-1}{m-1}z^n.$$

Comparing coefficients and simplifying shows that (1.28) is equivalent with (1.3).

Thus the polynomials $l_n(x,m,s)$ are multiples of convolutions of Fibonacci polynomials.

## 2. Generalized q-Fibonacci and q-Lucas polynomials

The polynomials $l_n(x,m,s)$ are a common generalization of Fibonacci, Lucas and Chebyshev polynomials, but in our search for $q-$analogues we were led to two quite different classes of polynomials (so as already in [6] for the case $m=0$ and $m=1$).

Let us start with the Fibonacci and Lucas case.

First some remarks about notation. As usual we set

$(a;q)_n = (1-a)(1-qa)\cdots(1-q^{n-1}a)$ with $(a;q)_0 = 1.$

Let $[n]=[n]_q = 1+q+\cdots+q^{n-1} = \frac{1-q^n}{1-q}$, $\begin{bmatrix}n\\k\end{bmatrix} = \frac{(q;q)_n}{(q;q)_k(q;q)_{n-k}}$ for $0\leq k\leq n$ and

$[n]! = \frac{(q;q)_n}{(1-q)^n}.$

The orthogonal polynomials $L_n(x,m,s,q)$ which have as moments the normalized $q-$ super Catalan numbers

$$\sigma_q(m,n) = \frac{[2n]!_q [m]!_q}{[n]!_q [m+n]!_q} \tag{2.1}$$

have ugly coefficients.

Therefore we consider instead the following polynomials whose coefficients are nice $q-$ analogues of the coefficients of $l_n(x,m,s)$.

$$l_n(x,m,s,q) = \sum_{k=0}^{\left\lfloor \frac{n}{2} \right\rfloor} s^k q^{\binom{k}{2}} \frac{[n]!_q}{[k]!_q [n-2k]!_q} \frac{[m+n-k-1]!_q}{[m+n-1]!_q} x^{n-2k}. \tag{2.2}$$

These polynomials are no longer orthogonal, but satisfy the following curious $q-$ analogue of the recursion (1.3).

$$l_n(x,m,s,q) = \left(x - (1-q)sq^{-m}D_q\right)l_{n-1}(x,m,s,q) + sq^{-m}\lambda_{n-2}(m,q)l_{n-2}(x,m,s,q) \tag{2.3}$$

with

$$\lambda_n(m,q) = \frac{[n+1]_q [n+2m]_q}{[n+m]_q [n+m+1]_q}. \tag{2.4}$$

Here $D_q$ is the $q-$ differentiation operator defined by $D_q f(x) = \frac{f(x) - f(qx)}{(1-q)x}$.

The simplest $q-$ analogues of $l_n(x,m,-1)$ are therefore the polynomials $l_n(x,m,-q^m,q)$ which satisfy

$$l_n(x,m,-q^m,q) = \left(x + (1-q)D_q\right)l_{n-1}(x,m,-q^m,q) - \lambda_{n-2}(m,q)l_{n-2}(x,m,-q^m,q). \tag{2.5}$$

Writing $l_n(x,m,-q^m,q) = \sum_{k=0}^{\left\lfloor \frac{n}{2} \right\rfloor} c(n,k) x^{n-2k}$ (2.5) is equivalent with

$$c(n,k) = c(n-1,k) + \left(1 - q^{n+1-2k}\right)c(n-1,k-1) - \lambda_{n-2}(m,q)c(n-2,k-1).$$

By setting $c(n,k) = (-1)^k q^{\binom{k}{2}+mk} \begin{bmatrix} n-k \\ k \end{bmatrix}_q \frac{[n]!_q}{[n-k]!_q} \frac{[m+n-k-1]!_q}{[m+n-1]!_q}$ this is easily verified by induction.

The first terms are

$$1,\ x,\ x^2 - \frac{q^m[2]}{[m+1]},\ x^3 - \frac{q^m[2][3]}{[m+2]}x,\ x^4 - \frac{q^m[3][4]}{[m+3]}x^2 + \frac{q^{2m+1}[3][4]}{[m+2][m+3]}, \cdots.$$

Define a linear functional $\Lambda_{m,q}$ on the polynomials by $\Lambda_{m,q}\left(l_n(x,m,-q^m,q)\right) = [n=0]$.

As above we get the representation

$$x^n = \sum_{k=0}^{\left\lfloor \frac{n}{2} \right\rfloor} (-s)^k \frac{[n]!_q}{[k]!_q[n-2k]!_q} \frac{[n+m-2k]!_q}{[n+m-k]!_q} l_{n-2k}(x,m,s,q) \tag{2.6}$$

which implies that $\Lambda_{m,q}\left(x^{2n+1}\right) = 0$ and

$$\Lambda_{m,q}\left(x^{2n}\right) = q^{mn} \frac{[2n]!_q [m]!_q}{[n]!_q [m+n]!_q} = q^{mn} \sigma_q(m,n). \tag{2.7}$$

For $m=0$ this is [6], (3.13) and for $m=1$ it is [6], (3.5).

For $m=1$ formula (2.6) is equivalent with

$$x^{2n} = \sum_{k=0}^{n} q^{n-k} \frac{[2n]![2k+1]}{[n-k]![n+1+k]!} l_{2k}(x,1,-q,q) = \sum_{k=0}^{n} \left( \begin{bmatrix} 2n \\ n-k \end{bmatrix} - \begin{bmatrix} 2n \\ n-k-1 \end{bmatrix} \right) l_{2k}(x,1,-q,q),$$

$$x^{2n-1} = \sum_{k=1}^{n} q^{n-k} \frac{[2n-1]![2k]}{[n-k]![n+k]!} l_{2k-1}(x,1,-q,q) = \sum_{k=1}^{n} \left( \begin{bmatrix} 2n-1 \\ n-k \end{bmatrix} - \begin{bmatrix} 2n-1 \\ n-k-1 \end{bmatrix} \right) l_{2k-1}(x,1,-q,q).$$

(2.8)

Generalizations in other directions are given in [7] and the literature cited there.

**Remarks**

For $|q|<1$ it is possible to consider the limit for $m \to \infty$ in (2.2). If we replace $s$ by $\dfrac{s}{1-q}$ we get the polynomials

$$h_n(x,s,q) = \sum_{k=0}^{\left\lfloor \frac{n}{2} \right\rfloor} s^k q^{\binom{k}{2}} c(n,k,q) x^{n-2k} \tag{2.9}$$

with $c(n,k,q) = \dfrac{[n]!_q}{[k]!_q [n-2k]!_q}$, which are $q-$analogues of a variant of the Hermite polynomials.

They satisfy the recurrence

$$h_n(x,s,q) = x h_{n-1}(x,s,q) + q^{n-2}(1+q)[n-1] s h_{n-2}(x,s,q) + (1-q) q^{n-3} s^2 [n-1][n-2][n-3] h_{n-4}(x,s,q).$$
$$(2.10)$$

This can be verified with the Mathematica implementation qZeil of the $q-$analogue of Zeilberger's algorithm by Axel Riese.

By (2.6) the inverse is given by

$$x^n = \sum_{k=0}^{\left\lfloor \frac{n}{2} \right\rfloor} (-s)^k c(n,k,q) h_{n-2k}(x,s,q).$$

This implies that the moments of the polynomials $h_n(x,-1,q)$ are $\Lambda(x^{2n}) = \dfrac{[2n]!}{[n]!}$.

The normalized $q-$super Catalan numbers $\sigma_q(m,n)$ can be written as the product

$$\sigma_q(m,n) = \dfrac{[2n]!}{[n]!} \cdot \dfrac{[m]!}{[m+n]!}.$$

For the second factor we have an analogous situation with $c(n,k,q) = \begin{bmatrix} \left\lfloor \dfrac{n}{2} \right\rfloor \\ k \end{bmatrix}$.

For the polynomials

$$R(n,m,q) = \sum_{k=0}^{\lfloor \frac{n}{2} \rfloor} (-1)^k q^{\binom{k}{2}} \begin{bmatrix} \lfloor \frac{n}{2} \rfloor \\ k \end{bmatrix} \frac{1}{\prod_{j=1}^{k}[m+n-j]} x^{n-2k}$$

we get

$$x^n = \sum_{k=0}^{\lfloor \frac{n}{2} \rfloor} \begin{bmatrix} \lfloor \frac{n}{2} \rfloor \\ k \end{bmatrix} \frac{1}{\prod_{j=0}^{k-1}[m+n-k-j]} R_{n-2k}(x,s,q).$$

Therefore their moments are $\dfrac{[m]!}{[m+n]!}$.

By letting $m \to \infty$ we get the polynomials

$$r_n(x,q) = \sum_{k=0}^{\lfloor \frac{n}{2} \rfloor} (-1)^k q^{\binom{k}{2}} \begin{bmatrix} \lfloor \frac{n}{2} \rfloor \\ k \end{bmatrix} x^{n-2k}$$

whose moments are 1.

Observe that $r_{2n}(x,q) = \prod_{j=0}^{n-1}(x^2 - q^j)$ and $r_{2n+1}(x,q) = x\prod_{j=0}^{n-1}(x^2 - q^j)$.

As is well known (cf. [3] and the literature cited there) formula (1.6) has an interesting $q-$analogue: Let $r(n) = \dfrac{n(3n-1)}{2}$ denote the pentagonal numbers. Then

$$(l_n(1,1,-1,q))_{n\geq 0} = (1,1,0,-q,-q^2,0,q^5,q^7,0,-q^{12},-q^{15},0,\cdots). \tag{2.11}$$

Thus

$$\begin{aligned} l_{3n}(1,1,-1,q) &= (-1)^n q^{r(n)} = (-1)^n q^{\frac{n(3n-1)}{2}}, \\ l_{3n+1}(1,1,-1,q) &= (-1)^n q^{r(-n)} = (-1)^n q^{\frac{n(3n+1)}{2}} \\ l_{3n+2}(1,1,-1,q) &= 0. \end{aligned} \tag{2.12}$$

qZeil gives the recurrences

$$l_{3n+i}(1,1,-1,q) + q^{3n-2+i} l_{3n-3+i}(1,1,-1,q) = 0$$

which imply (2.12).

Observing that $l_n(x,0,-1,q) = l_n(x,1,-1,q) - q^{n-1}l_{n-2}(x,1,-1,q)$ we see that

$$(l_n(1,0,-1,q))_{n\geq 0} = \left(1,1,-q,-q-q^2,-q^2,q^5,q^5+q^7,-q^{12},\cdots\right). \tag{2.13}$$

We also get nice formulae for $m=2$.

Using $q-$Zeil we see that for $i \in \{0,1,2\}$ $f(n,i) = l_{3n+i}(1,2,-1,q)[3n+i+1]$ satisfies the recurrence

$$f(n,i) + q^{3n-2+i}[3]f(n-1,i) + q^{6n-6+2i}[3]f(n-2,i) + q^{9n-12+3i}f(n-3,i) = 0.$$

Therefore

$$(l_n(1,2,-1,q)[n+1])_{n\geq 0} = \left(1, 1+q, q^2, -q(1+q), -q^2(1+q)(1+q^2), -q^6(1+q), q^5(1+q+q^2), \cdots\right).$$
(2.14)

This means

$$l_{3n}(1,2,-1,q)[3n+1] = (-1)^n q^{\frac{n(3n-1)}{2}} [n+1],$$

$$l_{3n+1}(1,2,-1,q)[3n+2] = (-1)^n q^{\frac{n(3n+1)}{2}} [n+1]\left(1+q^{n+1}\right),$$

$$l_{3n+2}(1,2,-1,q)[3n+3] = (-1)^n q^{\frac{3n^2+5n+4}{2}} [n+1].$$

## 3. Generalized q-Chebyshev polynomials

The most natural $q-$analogues of the Chebyshev-like polynomials $l_n\left(x,m,-\frac{1}{4}\right)$ turn out to be the polynomials $v_n(x,m,1,q)$ where

$$v_n(x,m,s,q) = \sum_{k=0}^{\left\lfloor \frac{n}{2} \right\rfloor} (-s)^k q^{k^2} \frac{[n]!}{[k]![n-2k]!} \frac{[m+n-k-1]!}{[m+n-1]!} \frac{1}{(-q;q)_k (-q^{n+m-k};q)_k} x^{n-2k}. \tag{3.1}$$

The first terms of the sequence $(v_n(x,m,s,q))_{n \geq 0}$ are

$$1, \; x, \; x^2 - \frac{qs}{[2m+2]}, \; x^3 - \frac{q[3]}{[2m+4]}sx, \; x^4 - \frac{q(1+q^2)[3]}{[2m+6]}sx^2 + \frac{q^4[3]s^2}{[2m+4][2m+6]}, \ldots.$$

They satisfy

$$v_n(x,m,s,q) = xv_{n-1}(x,m,s,q) - s\lambda_{n-2}(m,q)\frac{q^{n-1}}{(1+q^{n+m-2})(1+q^{n+m-1})}v_{n-2}(x,m,s,q) \quad (3.2)$$

with initial values

$v_0(x,m,s,q) = 1$ and $v_1(x,m,s,q) = x$.

For $m=0$ and $m=1$ they reduce to the monic $q-$Chebyshev polynomials (cf. [6]).

Since they satisfy the 3-term recurrence (3.2) they are orthogonal with respect to the linear functional $\Phi_{m,q}$ defined by $\Phi_{m,q}(v_n(x,m,1,q)) = [n=0]$.

We have

$$x^n = \sum_{k=0}^{\lfloor \frac{n}{2} \rfloor} \frac{[n]!}{[k]![n-2k]!} \frac{[n+m-2k]!}{[n+m-k]!} \frac{q^k s^k}{(-q;q)_k(-q^{n+m+1-2k};q)_k} v_{n-2k}(x,m,s,q). \quad (3.3)$$

This gives the moments

$$\Phi_{m,q}(x^{2n}) = \begin{bmatrix} 2n \\ n \end{bmatrix} \frac{[n]![m]!}{[n+m]!} \frac{q^n}{(-q;q)_n(-q^{m+1};q)_n} = \frac{q^n}{(-q;q)_n(-q^{m+1};q)_n}\sigma_q(m,n). \quad (3.4)$$

Also here we can consider the limit for $m \to \infty$ in (3.1) and get a version of the discrete $q-$ Hermite polynomials (cf. [2])

$$H_n(x,s,q) = \sum_{k=0}^{\lfloor \frac{n}{2} \rfloor} (-s)^k q^{k^2} \frac{[n]!}{[k]![n-2k]!} \frac{1}{(-q;q)_k} x^{n-2k}, \qquad (3.5)$$

whose moments are $\dfrac{(qs)^n}{(-q,q)_n} \dfrac{[2n]!}{[n]!} = (qs)^n [2n-1]!!.$

The sequence $\left( \dfrac{[m]!}{[m+n]!(-q^{m+1},q)_n} \right) = \left( \dfrac{1}{\prod_{j=1}^{n}[2m+2j]_q} \right)$ is the moment sequence of the polynomials

$$\sum_{k=0}^{\lfloor \frac{n}{2} \rfloor} (-1)^k q^{2\binom{k}{2}} \begin{bmatrix} \lfloor \frac{n}{2} \rfloor \\ k \end{bmatrix}_{q^2} \frac{1}{\prod_{j=1}^{k}[2(m+n-j)]_q}.$$

Let now $b(n,k,m,q)$ be the weight of the set of all non-negative lattice paths from $(0,0)$ to $(n,k)$ consisting of up-steps $(i,j) \to (i+1,j+1)$ with weight 1 and down-steps $(i,j+1) \to (i+1,j)$ with weight $\mu_j(m,q) = \lambda_j(m,q) \dfrac{q^{j+1}}{(1+q^{j+m})(1+q^{j+m+1})}$ and $c(\ell,k,m,q)$ be the weight of the set of all non-negative lattice paths from $(n,k)$ to $(n+\ell,0)$.

We know that

$$b(2n,2k,m,q) = a(2n,n-k,m,q) = \begin{bmatrix} 2n \\ n-k \end{bmatrix} \frac{[n+k]![m+2k]!}{[2k]![n+m+k]!} \frac{q^{n-k}}{(-q;q)_{n-k}(-q^{m+1+2k};q)_{n-k}}$$

It is easily verified that

$c(2n,0,m,q) = b(2n,0,m,q)$ and for $k > 0$

$$c(2n,2k,m,q) = b(2n,2k,m,q) \frac{q^{k+2k^2}(q;q)_{2k}(q^{2m+1};q)_{2k-1}}{(q^{m+1};q)_{2k-1}^2(-q^{m+1};q)_{2k-1}^2(1-q^{2m+4k})}$$

$$= \begin{bmatrix} 2n \\ n-k \end{bmatrix} \frac{[n+k]![m+2k]!}{[2k]![n+m+k]!} \frac{q^{k+2k^2}(q;q)_{2k}(q^{2m+1};q)_{2k-1}}{(q^{m+1};q)_{2k-1}^2(-q^{m+1};q)_{2k-1}^2(1-q^{2m+4k})}.$$

For example for $m=0$ we get

$$b(2n,2k,0,q) = \begin{bmatrix} 2n \\ n-k \end{bmatrix} \frac{q^{n-k}}{(-q;q)_{n-k}(-q^{1+2k};q)_{n-k}} \quad \text{and}$$

$$c(2n,2k,0,q) =$$

$$= \begin{bmatrix} 2n \\ n-k \end{bmatrix} \frac{q^{n+2k^2}}{(-q;q)_{n-k}(-q^{1+2k};q)_{n-k}(-q;q)_{2k-1}(-q;q)_{2k}}.$$

Thus we get as a $q-$analogue of (1.18)

$$\begin{bmatrix} 2n+2\ell \\ n+\ell \end{bmatrix} \frac{1}{(-q;q)_{\ell+n}^2} = \begin{bmatrix} 2n \\ n \end{bmatrix} \begin{bmatrix} 2\ell \\ \ell \end{bmatrix} \frac{1}{(-q;q)_n^2(-q;q)_\ell^2}$$
$$+ \sum_{k=1}^{\min(\ell,n)} \begin{bmatrix} 2n \\ n-k \end{bmatrix} \begin{bmatrix} 2\ell \\ \ell-k \end{bmatrix} \frac{q^{2k^2-k}(1+q^{2k})}{(-q;q)_{n-k}(-q;q)_{n+k}(-q;q)_{\ell-k}(-q;q)_{\ell+k}}. \tag{3.6}$$

As $q-$ analogue of (1.23) we get

$$\begin{bmatrix} 2n \\ n \end{bmatrix} \frac{1}{(-q;q)_n^2} + \sum_{k=1}^n (-1)^k q^{2\binom{k}{2}} \begin{bmatrix} 2n \\ n-k \end{bmatrix} \frac{1+q^{2k}}{(-q;q)_{n-k}(-q;q)_{n+k}} = [n=0]. \tag{3.7}$$

An analogue of (1.24) is

$$\sum_{k=0}^n (-1)^k q^{2\binom{k}{2}} \begin{bmatrix} 2n \\ n-k \end{bmatrix} \frac{[2k+1]}{[n+k+1]} \frac{1+q^{2k+1}}{(-q;q)_{n-k}(-q;q)_{n+k+1}} = [n=0]. \tag{3.8}$$

In general formula (1.21) gives $\sum_{k=0}^{n}(-1)^k b(2n,2k,m,q)\prod_{j=0}^{k-1}\lambda_{2j}(m,q)=[n=0]$.

Since

$$\prod_{j=0}^{k-1}\mu_{2j}(m,q)=s(n+m-1,n,q)\frac{q^{n^2}}{(-q;q)_n(-q^{n+m};q)_n}=\frac{[2n]![n+m-1]!}{[n]![2n+m-1]!}\frac{q^{n^2}}{(-q;q)_n(-q^{n+m};q)_n}$$

we get by simplifying that this identity is equivalent with

$$\sum_{k=0}^{n}(-1)^k\begin{bmatrix}n\\k\end{bmatrix}_{q^2}\frac{[m+2k]_{q^2}\,q^{2\binom{k}{2}}[k+m-1]_{q^2}!}{[n+m+k]_{q^2}!}=[n=0]. \quad (3.9)$$

As analogue of (1.26) we get

$$\sum_{k=0}^{n}b(2n,2k,m,q)\prod_{j=0}^{k-1}\mu_{2j}(m,q)=q^n\frac{[2n]_{q^2}!\,(-1;q^2)_n}{[n]_{q^2}!\,(-q;q)_{2n}}\frac{1}{[m+1]_{q^2}[m+3]_{q^2}\cdots[m+2n-1]_{q^2}}.$$
(3.10)

By simplifying

$$\sum_{k=0}^{n}\begin{bmatrix}2n\\n-k\end{bmatrix}\frac{[n+k]![m+2k]!}{[2k]![n+m+k]!}\frac{q^{n-k}}{(-q;q)_{n-k}(-q^{m+1+2k};q)_{n-k}}\frac{[2k]![k+m-1]!}{[k]![2k+m-1]!}\frac{q^{k^2}}{(-q;q)_k(-q^{k+m};q)_k}$$

$$=\frac{[2n]!}{[n]_2!}\sum_{k=0}^{n}\begin{bmatrix}n\\k\end{bmatrix}_2\frac{[m+2k]_{q^2}\,q^{n-k+k^2}[k+m-1]_{q^2}!}{[n+m+k]_2!}=q^n\frac{[2n]_{q^2}!\,(-1;q^2)_n}{[n]_{q^2}!\,(-q;q)_{2n}}\frac{1}{[m+1]_{q^2}[m+3]_{q^2}\cdots[m+2n-1]_{q^2}}$$

this identity is equivalent with

$$\sum_{k=0}^{n}\begin{bmatrix}n\\k\end{bmatrix}q^{\binom{k}{2}}[m+2k]\frac{[m+k-1]!}{[m+n+k]!}=\frac{(-1;q)_n}{[m+1][m+3]\cdots[m+2n-1]}. \quad (3.11)$$

The Mathematica implementation qZeil of the $q-$ analogue of Zeilberger's algorithm gives a simple computer proof:

```
qZeil[qBinomial[n, k, q] qBinomial[m + 2 k, 1, q] q^Binomial[k, 2]
  qFactorial[k + m - 1, q] / qFactorial[n + m + k, q], {k, 0, n}, n, 1, {m}]
```

$$\operatorname{SUM}[n] == -\frac{\left(1 - \frac{1}{q}\right) q \left(1 + q^{-1+n}\right) \operatorname{SUM}[-1+n]}{1 - q^{-1+m+2n}}$$

The following human proof uses an idea of [10].

Let $a = q^m$. Then (3.11) can be written as

$$S_n(a) = \sum_{k=0}^{n} \begin{bmatrix} n \\ k \end{bmatrix} q^{\binom{k}{2}} \frac{\left(1 - q^{2k} a\right)}{\left(q^k a; q\right)_{n+1}} = \frac{(-1; q)_n}{\left(qa; q^2\right)_n}. \qquad (3.12)$$

Let

$$F_{n,k}(a) = \begin{bmatrix} n \\ k \end{bmatrix} q^{\binom{k}{2}} \frac{\left(1 - q^{2k} a\right)}{\left(q^k a; q\right)_{n+1}}.$$

Then we have

$$F_{n,k}(a) = \frac{F_{n-1,k}(a)}{1 - q^n a} + q^{n-1} \frac{F_{n-1,k-1}(q^2 a)}{1 - q^n a}. \qquad (3.13)$$

To show this we use the identities $\begin{bmatrix} n \\ k \end{bmatrix} = q^k \begin{bmatrix} n-1 \\ k \end{bmatrix} + \begin{bmatrix} n-1 \\ k-1 \end{bmatrix}$ and $\begin{bmatrix} n \\ k \end{bmatrix} = \begin{bmatrix} n-1 \\ k \end{bmatrix} + q^{n-k} \begin{bmatrix} n-1 \\ k-1 \end{bmatrix}$

which can be combined to give

$$\begin{bmatrix} n \\ k \end{bmatrix} (1 - q^n a) = \begin{bmatrix} n-1 \\ k \end{bmatrix} (1 - q^{n+k} a) + \begin{bmatrix} n-1 \\ k-1 \end{bmatrix} (1 - q^k a) q^{n-k}.$$

This implies (3.13):

$$F_{n,k}(a) = \begin{bmatrix} n \\ k \end{bmatrix} q^{\binom{k}{2}} \frac{\left(1 - q^{2k} a\right)}{\left(q^k a; q\right)_{n+1}} = \frac{\begin{bmatrix} n-1 \\ k \end{bmatrix}(1 - q^{n+k} a)}{1 - q^n a} q^{\binom{k}{2}} \frac{\left(1 - q^{2k} a\right)}{\left(q^k a; q\right)_{n+1}} + \frac{\begin{bmatrix} n-1 \\ k-1 \end{bmatrix}(1 - q^k a) q^{n-k}}{1 - q^n a} q^{\binom{k}{2}} \frac{\left(1 - q^{2k} a\right)}{\left(q^k a; q\right)_{n+1}}$$

$$= \frac{F_{n-1,k}(a)}{1 - q^n a} + \frac{q^{n-1}}{1 - q^n a} \begin{bmatrix} n-1 \\ k-1 \end{bmatrix}(1 - q^k a) q^{\binom{k-1}{2}} \frac{\left(1 - q^{2k-2} q^2 a\right)}{\left(q^{k-1} q^2 a; q\right)_n} = \frac{F_{n-1,k}(a)}{1 - q^n a} + \frac{q^{n-1}}{1 - q^n a} F_{n-1,k-1}(q^2 a).$$

Therefore

$$S_n(a) = \frac{1}{1 - q^n a} S_{n-1}(a) + \frac{q^{n-1}}{1 - q^n a} S_{n-1}(q^2 a).$$

The right-hand side satisfies the same recurrence:

$$\frac{(-1;q)_n}{(qa;q^2)_n} = \frac{1}{1-q^n a} \frac{(-1;q)_{n-1}}{(qa;q^2)_{n-1}} + \frac{q^{n-1}}{1-q^n a} \frac{(-1;q)_{n-1}}{(q^3 a;q^2)_{n-1}}.$$

For this can be simplified to the trivial identity $(1+q^{n-1})(1-q^n a) = (1-q^{2n-1}a) + q^{n-1}(1-qa)$.

Since (3.12) is true for $n=0$ the proof is complete.

**4. The polynomials $v_n(x,m,s)$ and their q-analogues as polynomials in $s$.**

Consider the bivariate polynomials in $x$ and $s$

$$v_n(x,m,s) = \sum_{k=0}^{\lfloor \frac{n}{2} \rfloor} (-1)^k \frac{n!}{(k)!(n-2k)!} \frac{(m+n-k-1)!}{(m+n-1)!} \frac{1}{2^{2k}} s^k x^{n-2k} = l_n\left(x,m,-\frac{s}{4}\right) \quad (4.1)$$

The first terms are

$$1, x, x^2 - \frac{s}{2m+2}, x^3 - \frac{3sx}{2m+4}, x^4 - \frac{3sx^2}{m+3} + \frac{3s^2}{(2m+4)(2m+6)}, x^5 - \frac{5sx^3}{m+4} + \frac{15s^2 x}{(2m+6)(2m+8)}, \cdots.$$

Since both $(v_{2n}(x,m,s))_{n\geq 0}$ and $(v_{2n+1}(x,m,s))_{n\geq 0}$ are bases for the vector space of polynomials in $s$ it is clear that $v_{2n+1}(x,m,s)$ has an expansion of the form

$$v_{2n+1}(x,m,s) = \sum_{k=0}^{n} a(n,k,m) x^{2k+1} v_{2n-2k}(x,m,s). \quad (4.2)$$

More precisely we show that for $m > 0$

$$v_{2n+1}(x,m,s) = \sum_{k=0}^{n} (-1)^k \binom{2n+1}{2k+1} \frac{(m+2n-2k-1)!(2m+2n)!}{(m+2n)!(2m+2n-2k-1)!} \frac{T_{2k+1}}{2^{2k+1}} x^{2k+1} v_{2n-2k}(x,m,s),$$

(4.3)

where $T_{2n+1}$ is a tangent number defined by

$$\frac{e^z - e^{-z}}{e^z + e^{-z}} = \sum_{n\geq 0} (-1)^n \frac{T_{2n+1}}{(2n+1)!} z^{2n+1}. \quad (4.4)$$

The first terms are $(T_{2n+1})_{n\geq 0} = (1, 2, 16, 272, 7936, \cdots)$.

For $m = 0$ we have the following coefficients (cf. [5]):

$$a(n,n,0) = (-1)^n \frac{T_{2n+1}}{2^{2n}} \qquad (4.5)$$

and for $k < n$

$$a(n,k,0) = (-1)^k \binom{2n+1}{2k+1} \frac{T_{2k+1}}{2^{2k+1}}. \qquad (4.6)$$

To give a direct proof for this special case consider the polynomials $w_n(s)$ defined by $w_n(s) = v_n(1,0,s)$ for $n > 0$ and $w_0(s) = 2$. They satisfy $w_n(s) = w_{n-1}(s) - \frac{s}{4} w_{n-2}(s)$ with initial values $w_0(s) = 2$ and $w_1(s) = 1$. Binet's formula gives $w_n(s) = \alpha^n + \beta^n$ with 
$\alpha = \frac{1+\sqrt{1-s}}{2}$ and $\beta = \frac{1-\sqrt{1-s}}{2}$.

Therefore

$$V(s,z) := \sum_{n \geq 0} w_n(s) \frac{z^n}{n!} = \sum_{n \geq 0} (\alpha^n + \beta^n) \frac{z^n}{n!} = e^{\alpha z} + e^{\beta z} = e^{(1-\beta)z} + e^{(1-\alpha)z} = e^z \sum_{n \geq 0} (-1)^n (\alpha^n + \beta^n) \frac{z^n}{n!}$$
$$= e^z V(s,-z)$$

Thus

$$\sum_{n \geq 0} w_{2n+1}(s) \frac{z^{2n+1}}{(2n+1)!} = \frac{V(s,z) - V(s,-z)}{2} = \frac{V(s,z)}{2}(1 - e^{-z}) \text{ and}$$

$$\sum_{n \geq 0} w_{2n}(s) \frac{z^{2n}}{(2n)!} = \frac{V(s,z) + V(s,-z)}{2} = \frac{V(s,z)}{2}(1 + e^{-z}) \text{ and therefore}$$

$$\sum_{n \geq 0} w_{2n+1}(s) \frac{z^{2n+1}}{(2n+1)!} = \frac{1-e^{-z}}{1+e^{-z}} \sum_{n \geq 0} w_{2n}(s) \frac{z^{2n}}{(2n)!} = \sum_{n \geq 0} (-1)^n \frac{T_{2n+1}}{(2n+1)!} \left(\frac{z}{2}\right)^{2n+1} \sum_{n \geq 0} w_{2n}(s) \frac{z^{2n}}{(2n)!}.$$

Comparing coefficients gives

$$w_{2n+1}(s) = \sum_{k=0}^{n} \binom{2n+1}{2k+1} \frac{T_{2k+1}}{2^{2k+1}} w_{2n-2k}(s)$$

and thus (4.5) and (4.6).

For example

$$v_5(1,0,s) = \binom{5}{1}\frac{T_1}{2}v_4(1,0,s) - \binom{5}{3}\frac{T_3}{8}v_2(1,0,s) + \binom{5}{5}\frac{T_5}{16} = \frac{5}{2}v_4(1,0,s) - \frac{5}{2}v_2(1,0,s) + 1.$$

In general we claim that for $m > 0$

$$a(n,k,m) = (-1)^k \binom{2n+1}{2k+1}\frac{(m+2n-2k-1)!(2m+2n)!}{(m+2n)!(2m+2n-2k-1)!}\frac{T_{2k+1}}{2^{2k+1}}. \qquad (4.7)$$

By comparing coefficients this is equivalent with

$$\sum_{n\geq 0}\frac{(m+2n)!}{(2m+2n)!}v_{2n+1}(1,m,s)\frac{z^{2n+1}}{(2n+1)!}$$
$$= \sum_{k\geq 0}(-1)^k\frac{T_{2k+1}}{2^{2k+1}}\frac{z^{2k+1}}{(2k+1)!}\sum_{n\geq 0}\frac{(m+2\ell-1)!}{(2m+2\ell-1)!}v_{2\ell}(1,m,s)\frac{z^{2\ell}}{(2\ell)!}. \qquad (4.8)$$

Let now

$$V_m(z) = V_m(z,s) = \sum_{n\geq 0}\frac{(m+n-1)!}{(2m+n-1)!}v_n(1,m,s)\frac{z^n}{n!}. \qquad (4.9)$$

Then (4.8) can be written as

$$V_m(z) - V_m(-z) = \frac{1-e^{-z}}{1+e^{-z}}\left(V_m(z) - V_m(-z)\right) \qquad (4.10)$$

or

$$V_m(-z) = e^{-z}V_m(z). \qquad (4.11)$$

This is equivalent with

$$\sum_{k=0}^{n}(-1)^{k}\binom{n}{k}\frac{(m+k-1)!}{(2m+k-1)!}v_{k}(1,m,s) = \frac{(m+n-1)!}{(2m+n-1)!}v_{n}(1,m,s). \tag{4.12}$$

Comparing coefficients of $s^j$ this is again equivalent with

$$\sum_{k=2j}^{n}(-1)^{k}\frac{(2m+n-1)!(n-2j)!(m+k-j-1)!}{(n-k)!(2m+k-1)!(k-2j)!(m+n-j-1)!} = 1 \tag{4.13}$$

or with

$$\sum_{k=0}^{n}(-1)^{k}\frac{(m+k+j-1)!}{(n-k)!(2m+k+2j-1)!(k)!} = \frac{(m+n+j-1)!}{(2m+n+2j-1)!(n)!}.$$

The left-hand side can be written in the form

$$\frac{(m+j-1)!}{n!(2m+2j-1)!}\sum_{k=0}^{n}\frac{(m+j)_{k}(-n)_{k}}{(2m+2j)_{k}k!}$$

by setting $(a)_{n} = a(a+1)\cdots(a+n-1)$.

By Gauß's theorem on hypergeometric series we have

$$\sum_{k=0}^{n}\frac{(a)_{k}(b)_{k}}{(c)_{k}k!} = \frac{(c-1)!(c-a-b-1)!}{(c-a-1)!(c-b-1)!}.$$

Therefore we get

$$\sum_{k=0}^{n}(-1)^{k}\frac{(m+k+j-1)!}{(n-k)!(2m+k+2j-1)!(k)!} = \frac{(m+j-1)!}{n!(2m+2j-1)!}\sum_{k=0}^{n}\frac{(m+j)_{k}(-n)_{k}}{(2m+2j)_{k}k!}$$
$$= \frac{(m+j-1)!}{n!(2m+2j-1)!}\frac{(2m+2j-1)!(m+j+n-1)!}{(m+j-1)!(2m+2j+n-1)!} = \frac{(m+j+n-1)!}{n!(2m+2j+n-1)!},$$

which gives (4.13).

We can also give a computer proof with fastZeil:

```
Zb[(-1)^k (2 m + n - 1)! (n - 2 j)! (m + k - j - 1)!/((n - k)! (2 m + k - 1)! (k - 2 j)! (m + n - j - 1)!),
 {k, 2 j, n}, n, 1]
```

If `-2 j + n' is a natural number, then:

{SUM[n] − SUM[1 + n] == 0}

Since SUM[0] = 1 the assertion follows.

For $m = 1$ we have

$$a(n,k,1) = (-1)^k \binom{2n+2}{2k+1} \frac{T_{2k+1}}{2^{2k+1}} = (-1)^k \binom{2n+2}{2k+2} \frac{1}{2n-2k+1} G_{2k+2}. \quad (4.14)$$

The Genocchi numbers $(G_{2n})_{n \geq 0} = (0, 1, 1, 3, 17, 155, 2073, \cdots)$ are given by their generating function

$$z \frac{e^z - e^{-z}}{e^z + e^{-z}} = \sum_{n \geq 0} (-1)^{n-1} 2^{2n-1} \frac{G_{2n}}{(2n)!} z^{2n}. \quad (4.15)$$

Note that

$$G_{2n+2} = \frac{(n+1)T_{2n+1}}{2^{2n}}.$$

For example

$$v_5(1,1,s) = \binom{6}{2} \frac{G_2}{5} v_4(1,1,s) - \binom{6}{4} \frac{G_4}{3} v_2(1,1,s) + \binom{6}{6} \frac{G_6}{1}$$

$$= \binom{6}{1} \frac{T_1}{2} v_4(1,1,s) - \binom{6}{3} \frac{T_3}{2^3} v_2(1,1,s) + \binom{6}{5} \frac{T_5}{2^5} = 3v_4(1,1,s) - 5v_2(1,1,s) + 3.$$

**Remark**

Identity (4.11) holds for all $s$, therefore especially for $s = 0$. By (4.9)

$$V_m(z,0) = \sum_{n \geq 0} \frac{(m+n-1)!}{(2m+n-1)!} \frac{z^n}{n!} = \frac{m!}{(2m)!} \sum_{n \geq 0} s(n+m-1,m) \frac{z^n}{n!}.$$

The simplest functions satisfying $e^z f(-z) = f(z)$ are $z^{2k}(e^z + 1)$ and $z^{2k+1}(e^z - 1)$.

Comparing coefficients shows that $V_m(z,0)$ is a linear combination of such functions.

We get

$$V_1(z,0) = \frac{e^z - 1}{z},$$

$$V_2(z,0) = \frac{z(e^z + 1) - 2(e^z - 1)}{z^3},$$

$$V_3(z,0) = \frac{z^2(e^z - 1) - 6z(e^z + 1) + 12(e^z - 1)}{z^5},$$

and generally

$$V_m(z,0) = \frac{1}{z^{2m-1}} \sum_{k=0}^{m-1} (-1)^k \binom{m+k-1}{2k} 2^k (2k-1)!! z^{m-1-k} \left(e^z + (-1)^{m+k}\right). \tag{4.16}$$

**The q-analogues**

Now consider the $q$-analogues

$$v_n(x,m,s,q) = \sum_{k=0}^{\lfloor \frac{n}{2} \rfloor} (-1)^k q^{k^2} \begin{bmatrix} n-k \\ k \end{bmatrix} \frac{[n]!}{[n-k]!} \frac{[m+n-k-1]!}{[m+n-1]!} \frac{1}{(-q;q)_k \left(-q^{n+m-k};q\right)_k} s^k x^{n-2k}. \tag{4.17}$$

They satisfy the recursion

$$v_n(x,m,s,q) = xv_{n-1}(x,m,s,q) - s\lambda_{n-2}(m,q) \frac{q^{n-1}}{\left(1+q^{n+m-2}\right)\left(1+q^{n+m-1}\right)} v_{n-2}(x,m,s,q).$$

**Remark**

A $q$-analogue of (1.8) is

$$v_n\left(1,m,\frac{1}{q},q\right) = q^{\binom{n}{2}} \frac{[2m+n-1]}{[2m+2n-2]} \frac{[2m+n-2]}{[2m+2n-4]} \frac{[2m+n-3]}{[2m+2n-6]} \cdots \frac{[2m+1]}{[2m+2]}. \tag{4.18}$$

This can be proved by using the recurrence relation or the q-version of Zeilberger's algorithm.

As analogue of (4.3) for $m>0$ is

$$v_{2n+1}(x,m,s,q) = \sum_{k=0}^{n} (-1)^k \begin{bmatrix} 2n+1 \\ 2k+1 \end{bmatrix} \frac{[2m+2n]![m+2n-2k-1]!}{[m+2n]![2m+2n-2k-1]!} \frac{T_{2k+1}(q)x^{2k+1}}{\left(-q^{m+2n-2k};q\right)_{2k+1}} v_{2n-2k}(x,m,s,q). \tag{4.19}$$

The polynomials $T_{2k+1}(q)$ are the $q-$tangent numbers defined by

$$\sum_{n\geq 0}(-1)^n \frac{T_{2n+1}(q)}{[2n+1]!}z^{2n+1} = \frac{e_q(z)-e_q(-z)}{e_q(z)+e_q(-z)}, \quad (4.20)$$

where $e_q(z) = \sum_{n\geq 0}\frac{z^n}{[n]!}$ is the $q-$exponential series.

The first terms of the sequence $T_{2n+1}(q)$ are

$$1, q(1+q), q^2(1+q)^2(1+q^2)^2, q^3(1+q)^2(1+q^2)(1+q^3)(1+q+3q^2+2q^3+3q^4+2q^5+3q^6+q^7+q^8), \cdots$$

Let

$$V_m(z,q) = V_m(z,q,s) = \sum_{n\geq 0}\frac{[m+n-1]!(-q^m;q)_n}{[2m+n-1]!}v_n(1,m,s,q)\frac{z^n}{[n]!}. \quad (4.21)$$

Then (4.19) is equivalent with

$$V_m(z)-V_m(-z) = \frac{e_q(z)-e_q(-z)}{e_q(z)+e_q(-z)}(V_m(z)+V_m(-z)). \quad (4.22)$$

and thus with

$$\frac{e_q(z)}{e_q(-z)}V_m(-z,q) = V_m(z,q). \quad (4.23)$$

Since

$$\frac{e_q(z)}{e_q(-z)} = \sum_{n\geq 0}(-1;q)_n\frac{z^n}{[n]!} \text{ this is equivalent with}$$

$$\sum_{k=0}^n(-1)^k\begin{bmatrix}n\\k\end{bmatrix}(-1;q)_{n-k}(-q^m;q)_k\frac{[m+k-1]!}{[2m+k-1]!}v_k(1,m,s,q) = \frac{[m+n-1]!(-q^m;q)_n}{[2m+n-1]!}v_n(1,m,s,q)$$
(4.24)

By comparing coefficients this again is equivalent with

$$\sum_{k=2j}^{n}(-1)^k \frac{(-1;q)_{n-k}\left(-q^m;q\right)_k\left(-q^{n+m-j};q\right)_j}{\left(-q^m;q\right)_n\left(-q^{k+m-j};q\right)_j} \frac{[2m+n-1]![m+k-j-1]![n-2j]!}{[n-k]![2m+k-1]![k-2j]![m+n-j-1]!}=1.$$

(4.25)

I want to thank Michael Schlosser for showing me that (4.25) follows from the $q-$Pfaff-Saalschütz formula (cf. [8], (II.12))

$$\sum_{k=0}^{n}\frac{(a;q)_k(b;q)_k\left(q^{-n};q\right)_k}{(c;q)_k\left(abc^{-1}q^{1-n};q\right)_k(q;q)_k}q^k = \frac{\left(\frac{c}{a};q\right)_n\left(\frac{c}{b};q\right)_n}{(c;q)_n\left(\frac{c}{ab};q\right)_n}.$$

(4.26)

First change $k$ to $k+2j$ and write all terms of the form $[p]!$ by $\frac{(q;q)_p}{(1-q)^p}$. Since

$(2m+n-1)+(m+k-j-1)+(n-2j)=(n-k)+(2m+k-1)+(k-2j)+(m+n-j-1)$

the powers of $1-q$ cancel and we get

$$\sum_{k=0}^{n-2j}(-1)^k \frac{(-1;q)_{n-k-2j}\left(-q^m;q\right)_{k+2j}\left(-q^{n+m-j};q\right)_j}{\left(-q^m;q\right)_n\left(-q^{k+m+j};q\right)_j} \frac{(q;q)_{2m+n-1}(q;q)_{m+k+j-1}(q;q)_{n-2j}}{(q;q)_{n-k-2j}(q;q)_{2m+k+2j-1}(q;q)_k(q;q)_{m+n-j-1}}$$

$$=\frac{\left(-q^m;q\right)_{2j}\left(-q^{n+m-j};q\right)_j(q;q)_{2m+n-1}(q;q)_{m+j-1}(q;q)_{n-2j}}{\left(-q^m;q\right)_n(q;q)_{2m+2j-1}(q;q)_{m+n-j-1}}$$

$$\times \sum_{k=0}^{n-2j}(-1)^k \frac{(-1;q)_{n-k-2j}}{(q;q)_{n-k-2j}} \frac{1}{\left(-q^{k+m+j};q\right)_j} \frac{\left(-q^{m+2j};q\right)_k\left(q^{m+j};q\right)_k}{\left(q^{2m+2j};q\right)_k(q;q)_k}$$

Now we apply the elementary identities ([8], (I.11) and (I.18))

$$\frac{(a;q)_{r-\ell}}{(b;q)_{r-\ell}}=\frac{(a;q)_r}{(b;q)_r}\frac{\left(\frac{q^{1-r}}{b};q\right)_\ell}{\left(\frac{q^{1-r}}{a};q\right)_\ell}\left(\frac{b}{a}\right)^\ell \quad \text{and} \quad \left(cq^k;q\right)_j=\frac{(c;q)_j\left(cq^j;q\right)_k}{(c;q)_k}$$

for $a=-1, b=q, r=n-2j, \ell=k, c=-q^{m+j}$ and obtain

$$\frac{\left(-q^{m};q\right)_{2j}\left(-q^{n+m-j};q\right)_{j}(q;q)_{2m+n-1}(q;q)_{m+j-1}(q;q)_{n-2j}}{\left(-q^{m};q\right)_{n}(q;q)_{2m+2j-1}(q;q)_{m+n-j-1}}$$

$$\times \sum_{k=0}^{n-2j} \frac{(-1;q)_{n-2j}\left(q^{2j-n};q\right)_{k}}{(q;q)_{n-2j}\left(-q^{1+2j-n};q\right)_{k}} q^{k} \frac{\left(-q^{m+j};q\right)_{k}}{\left(-q^{m+j};q\right)_{j}\left(-q^{m+2j};q\right)_{k}} \frac{\left(-q^{m+2j};q\right)_{k}\left(q^{m+j};q\right)_{k}}{\left(q^{2m+2j};q\right)_{k}(q;q)_{k}}$$

$$= \frac{\left(-q^{m};q\right)_{2j}\left(-q^{n+m-j};q\right)_{j}(-1;q)_{n-2j}}{\left(-q^{m};q\right)_{n}\left(-q^{m+j};q\right)_{j}} \frac{(q;q)_{2m+n-1}(q;q)_{m+j-1}}{(q;q)_{2m+2j-1}(q;q)_{m+n-j-1}}$$

$$\times \sum_{k=0}^{n-2j} \frac{\left(q^{m+j};q\right)_{k}\left(-q^{m+j};q\right)_{k}\left(q^{2j-n};q\right)_{k}}{\left(q^{2m+2j};q\right)_{k}\left(-q^{1+2j-n};q\right)_{k}(q;q)_{k}} q^{k}$$

By (4.26) we get

$$\sum_{k=0}^{n-2j} \frac{\left(q^{m+j};q\right)_{k}\left(-q^{m+j};q\right)_{k}\left(q^{2j-n};q\right)_{k}}{\left(q^{2m+2j};q\right)_{k}\left(-q^{1+2j-n};q\right)_{k}(q;q)_{k}} q^{k} = \frac{\left(q^{m+j};q\right)_{n-2j}\left(-q^{m+j};q\right)_{n-2j}}{\left(q^{2m+2j};q\right)_{n-2j}(-1;q)_{n-2j}}.$$

This gives (4.25) because our sum reduces to

$$\frac{\left(-q^{m};q\right)_{2j}\left(-q^{n+m-j};q\right)_{j}(-1;q)_{n-2j}}{\left(-q^{m};q\right)_{n}\left(-q^{m+j};q\right)_{j}} \frac{(q;q)_{2m+n-1}(q;q)_{m+j-1}}{(q;q)_{2m+2j-1}(q;q)_{m+n-j-1}} \frac{\left(q^{m+j};q\right)_{n-2j}\left(-q^{m+j};q\right)_{n-2j}}{\left(q^{2m+2j};q\right)_{n-2j}(-1;q)_{n-2j}}$$

$$= \frac{\left(-q^{m};q\right)_{2j}\left(-q^{n+m-j};q\right)_{j}\left(-q^{m+j};q\right)_{n-2j}}{\left(-q^{m};q\right)_{n}\left(-q^{m+j};q\right)_{j}} \frac{(q;q)_{2m+n-1}(q;q)_{m+j-1}\left(q^{m+j};q\right)_{n-2j}}{(q;q)_{2m+2j-1}(q;q)_{m+n-j-1}\left(q^{2m+2j};q\right)_{n-2j}} = 1 \cdot 1 = 1.$$

The Mathematica implementation qZeil also gives a simple computer proof of (4.25):

```
qZeil[(-1)^k qPochhammer[-1, q, n - k] qPochhammer[-q^m, q, k]
  qPochhammer[-q^(m + n - j), q, j] / qPochhammer[-q^(k + m - j), q, j] / qPochhammer[-q^m, q, n]
  qFactorial[2 m + n - 1, q] qFactorial[m + k - j - 1, q]
  qFactorial[n - 2 j, q] / qFactorial[n - k, q] / qFactorial[2 m + k - 1, q] / qFactorial[k - 2 j, q] /
   qFactorial[m + n - j - 1, q], {k, 0, n}, n, 1, {m, j}]

{SUM[n] == 1, {2 j - n ≠ 0, -1 + 2 m + n ≠ 0}}
```

We can also compute the series $V_{m}(z,q,0)$. If we set $E(z) = \dfrac{e_{q}(z)}{e_{q}(-z)}$ then we get

$$V_1(z,q,0) = \frac{E(z)-1}{2z},$$

$$V_2(z,q,0) = \frac{z(E(z)+1)-(E(z)-1)}{2q(1+q)z^3},$$

$$V_3(z,q,0) = \frac{qz^2(E(z)-1)-[3]z(E(z)+1)+[3](E(z)-1)}{2q^4(1+q)(1+q^2)z^5}$$

and in the general case

$$V_m(z,q,0) = \frac{\sum_{k=0}^{m-1}(-1)^k q^{\binom{m-1-k}{2}}\begin{bmatrix}m+k-1\\2k\end{bmatrix}[2k-1]!!z^{m-1-k}\left(E(z)+(-1)^{m+k}\right)}{q^{(m-1)^2}(-1;q)_m z^{2m-1}}. \tag{4.27}$$

**Remark**

Some of the above results give for $m \to \infty$ known results about discrete $q-$Hermite polynomials. Thus (4.21) reduces to the generating function

$$h(z,s,q) = \sum_{n \geq 0} H_n(1,s,q)\frac{z^n}{[n]!}.$$

In this case we have (cf. e.g. [2])

$$h(z,s,q) = \frac{e_q(z)}{e_{q^2}\left(\frac{qsz^2}{[2]}\right)}.$$

This implies (4.23) because $\dfrac{e_q(z)}{e_q(-z)}h(-z,s,q) = \dfrac{e_q(z)}{e_q(-z)}\dfrac{e_q(-z)}{e_{q^2}\left(\frac{qsz^2}{[2]}\right)} = \dfrac{e_q(z)}{e_{q^2}\left(\frac{qsz^2}{[2]}\right)} = h(z,s,q).$

Thus (4.19) gives

$$H_{2n+1}(x,s,q) = \sum_{k=0}^{n}(-1)^k\begin{bmatrix}2n+1\\2k+1\end{bmatrix}T_{2k+1}(q)x^{2k+1}H_{2n-2k}(x,s,q). \tag{4.28}$$